\documentclass[12pt,a4paper]{article}

\usepackage{amsfonts,amssymb,amstext,verbatim,,cancel,amsmath,amsthm,amsopn,amsfonts}
\usepackage{hyperref}

\usepackage{times}
\newtheorem{thm}{Theorem}

\newtheorem{lem}[thm]{Lemma}
\newtheorem{prop}[thm]{Proposition}
\newtheorem{proper}[thm]{Property}
\theoremstyle{definition}

\theoremstyle{remark}
\newtheorem{rem}[thm]{Remark}
\numberwithin{equation}{section}
\numberwithin{thm}{section}

\begin{document}

\title{Codes arising from incidence matrices of points and hyperplanes in $PG(n,q)$\footnote{The research was supported by Ministry for Education, University and Research of Italy MIUR (Project PRIN 2012 ``Geometrie di Galois e strutture di incidenza'') and by Italian National Group for Algebraic and Geometric Structures and their Applications (GNSAGA-INdAM).}}

\author{Olga Polverino, Ferdinando Zullo}

\maketitle

\begin{center}
  Seconda Universit\`a Degli Studi Di Napoli\\
  Dipartimento di Matematica e Fisica\\
  Via Vivaldi, 81100 Caserta\\
  Italy\\
  olga.polverino@unina2.it\\
  ferdinando.zullo@unina2.it
\end{center}

\section*{Abstract}

In this paper we completely characterize the words with second minimum weight in the $p-$ary linear code generated by the rows of the incidence matrix of points and hyperplanes of $PG(n,q)$, with $q=p^h$
and $p$ prime, proving that they are the scalar multiples of the difference of the incidence vectors of two distinct hyperplanes of $PG(n,q)$.

\section{Introduction}

Consider the projective space $PG(n,q)$, with $q=p^h$, $h \geq 1$ and $p$ prime.
We define the incidence matrix $A=(a_{i,j})$ of points and hyperplanes in $PG(n,q)$ as the matrix whose rows are indexed by the $\theta_n$ hyperplanes of $PG(n,q)$
and whose columns are indexed by the $\theta_n$ points of $PG(n,q)$, and with entry
$$ a_{i,j}= \left\{ \begin{array}{llrr} 1 & \text{if} \hspace{0.2cm} P_j \in H_i, \\ 0 & \text{otherwise}  \end{array} \right.. $$
The $p-$ary linear code of points and hyperplanes of $PG(n,q)$, which we denote by $ \mathcal{C}(n,q)$, is the code generated over
$\mathbb{F}_p$ by the rows of the matrix $A$.
These codes belong to a more general class of codes, the Reed-Muller codes.
For comprehensive references see e.g. \cite{As}, \cite{Giu} and \cite{Mc}.

\smallskip

The interest of these codes was born after the works of E. Prange \cite{Prange} and L. D. Rudolph \cite{Rudolph}, which showed that
projective planes can be used to define error-correcting codes.

\smallskip

Let $c$ be a codeword of $\mathcal{C}(n,q)$, the subset of $\{1,\ldots,\theta_n\}$ which corresponds to nonzero components of $c$ is said to be
the \emph{support of} $c$, and it will be denoted by $supp(c)$. We identify this set with the corresponding set of points of $PG(n,q)$.
The size of $supp(c)$ is said to be the \emph{weight} of $c$ and we will denote it by $wt(c)$.
Let $X$ be a subset of $PG(n,q)$, with $v^X$ we will denote the incidence vector of $X$.
Note that $wt(v^X)=|X|$.
Moreover, let $c_1=(a_1,\ldots,a_{\theta_n}),$ $c_2=(b_1,\ldots,b_{\theta_n}) \in \mathbb{F}_p^{\theta_n}$, the standard inner product is $\displaystyle (c_1,c_2)=\sum_{i=1}^{\theta_n} a_ib_i$. The orthogonal code is denoted by $\mathcal{C}(n,q)^\perp$ and is given by
$$ \mathcal{C}(n,q)^\perp=\{v \in \mathbb{F}_p^{\theta_n} : (v,c)=0 \hspace{0.1cm} \hspace{0.1cm} \forall \hspace{0.1cm} c \in \mathcal{C}(n,q) \}. $$
The \emph{hull} of $\mathcal{C}(n,q)$ is defined as $\mathcal{C}(n,q)\cap \mathcal{C}(n,q)^\perp$.

\smallskip

The fundamental parameters of these codes are well known (cfr. \cite{As} and \cite{Mc}):
the length is equal to the number of points of $PG(n,q)$, that is $\theta_n$;
the dimension is the p-rank of $A$, that is $\binom{p+n-1}{n}^h+1$;
the minimum distance is the number of points of a hyperplane, that is $\theta_{n-1}$.
Also,

\smallskip

\begin{thm}\label{mw}\cite{As}\cite{Mc}
The codewords of $\mathcal{C}(n,q)$ with minimum weight are the scalar multiples of the incidence vectors of hyperplanes.
\end{thm}

\smallskip

In \cite{C} the following result has been proved.

\smallskip

\begin{thm}\cite{C}\label{ch}
\begin{enumerate}
  \item In the $p-$ary linear code arising from $PG(2,p)$, $p$ prime, there are no codewords with weight in $]p+1,2p[$.
  \item The codewords of weight $2p$ in the $p-$ary linear code arising from $PG(2,p)$, $p$ prime, are the scalar multiples of the differences
        of the incidence vectors of two distinct lines of $PG(2,p)$.
\end{enumerate}
\end{thm}

\smallskip

In \cite{Lav1} the authors generalize 1. of Theorem \ref{ch} to codes generated by the rows of the incidence matrix of points and hyperplanes in $PG(n,q)$.
In particular, they obtain the following result.

\smallskip

\begin{thm}\cite[Corollary 19]{Lav1}\label{c19}
There are no codewords with weight in the open interval $]\theta_{n-1},2q^{n-1}[$ in the code $\mathcal{C}(n,q)$, $q=p^h$, $p$ prime, $p > 5$.
\end{thm}

\smallskip

In \cite{Zss}, the authors characterized small weight planar codewords of $\mathcal{C}(2,q)$ improving Theorem \ref{c19}, but no proof has been published yet.

\smallskip

In this paper we extend Results 1. and 2. of Chouinard (Theorem \ref{ch}) in $\mathcal{C}(n,q)$ for each prime power $q$.
More precisely, in Section \ref{new} we prove the following.

\smallskip

\begin{thm}\label{us}
Let $q=p^h$ with $p$ prime.
\begin{enumerate}
  \item There are no codewords of $\mathcal{C}(n,q)$ with weight in the interval $]\theta_{n-1},2q^{n-1}[$.
  \item The codewords of weight $2q^{n-1}$ in $\mathcal{C}(n,q)$ are the scalar multiples of the differences
        of the incidence vectors of two distinct hyperplanes of $PG(n,q)$.
\end{enumerate}
\end{thm}

\section{Preliminaries}

\subsection{Blocking sets}

Let $p$ a prime and $q=p^h$, with $h$ a positive integer.
A subset $B$ of $PG(n,q)$ is a $k-$\emph{blocking set} (or \emph{blocking set with respect to} $(n-k)$-\emph{subspaces}) of $PG(n,q)$, with
$1\leq k \leq n-1$, if each $(n-k)-$subspace intersects $B$ in at least one point.
If $k=1$, we simply say that $B$ is a blocking set of $PG(n,q)$. A $k-$blocking set is called \emph{trivial} if it contains a $k-$subspace.
An $(n-k)-$subspace which contains exactly one point of the $k-$blocking set $B$ is called $(n-k)-$\emph{tangent space} of $B$ and such a point
is called \emph{essential point}.
We say $B$ \emph{minimal} if each point of $B$ is an essential point for $B$.

\smallskip

An $(n-1)-$blocking set of $PG(n,q)$ small enough can be reduced in a unique way to a minimal blocking set, as proved in \cite{Lav}
by using Lemma 2.11 of \cite{Ga}. More precisely,

\smallskip

\begin{thm}\label{min}\cite[Corollary 1]{Lav}
Every $(n-1)-$blocking set in $PG(n,q)$, of size smaller than $q^{n-1}+\theta_{n-1}$, can be uniquely reduced to a minimal $(n-1)-$blocking set.
\end{thm}

\smallskip

\smallskip

\subsection{Blocking sets and codewords of small weight in $\mathcal{C}(n,q)$}

The following properties of the code $\mathcal{C}(n,q)$ can be easily verified.

\begin{proper}\cite[Lemmas 1, 2 and 3]{Lav}\label{3.3.5}\label{3.3.4}\label{3.3.3}
\begin{enumerate}
  \item If $U_1$ and $U_2$ are subspaces of dimension at least $1$ in $PG(n,q)$, then $v^{U_1}-v^{U_2} \in \mathcal{C}(n,q)^{\perp}$.
  \item The scalar product $(c,v^U)$, with $c \in \mathcal{C}(n,q)$ and $U$ an arbitrary subspace of dimension at least $1$, is a constant.
  \item A codeword $c$ is in $\mathcal{C}(n,q) \cap \mathcal{C}(n,q)^{\perp}$ if and only if $(c,v^U)=0$ for all
        subspaces $U$ with $dim(U)\geq1$.
\end{enumerate}
\end{proper}

\smallskip

Codewords of small weight in $\mathcal{C}(n,q)$ are related to $(n-1)-$blocking sets. Indeed in \cite{Lav}, generalizing Lemma 23 of \cite{C},
the authors prove the following.

\smallskip

\begin{thm}\label{3.3.11}\cite[Lemma 6]{Lav}
If $c \in \mathcal{C}(n,q)$, $c \neq \mathbf{0}$, with weight less than $2q^{n-1}$, then
\begin{enumerate}
  \item $supp(c)$ is a minimal $(n-1)-$blocking set of $PG(n,q)$;
  \item $c$, up to scalar, is the incidence vector of its support;
  \item $supp(c)$ intersects every line in $1 \pmod p$ points.
\end{enumerate}
\end{thm}

\smallskip

The next Theorem, due to A. Blokhuis, A. E. Brouwer and H. Wilbrink in \cite{BBH}, gives geometric information on codewords of $\mathcal{C}(2,q)$ with
components $0$ and $1$.

\begin{thm}\label{3.3.15}\cite[Proposition]{BBH}\label{BBW}
Let $X$ be a subset of $PG(2,q)$ such that $v^X \in \mathcal{C}(2,q)$ and let $Q$ be a point of $PG(2,q)$ such that $Q \notin X$. Then the points $P \in X$
for which the line $PQ$ is tangent to $X$ are collinear.
\end{thm}

\section{The second minimum weight of $\mathcal{C}(n,q)$ and the characterisation of the codeword of weight $2q^{n-1}$}\label{new}

In this section we prove Theorem \ref{us}.

\medskip

\begin{rem}\label{rm}\cite[Proof of Theorem 5]{Lav}
Note that the restriction of a codeword to a subspace is a codeword of the code associated with the fixed subspace.
Indeed, if $c \in \mathcal{C}(n,q)$ then there exist $\alpha_1,\ldots,\alpha_{\theta_n} \in \mathbb{F}_p$ such that
$c=\alpha_1v^{H_1}+\cdots+\alpha_{\theta_n}v^{H_{\theta_n}}$, where $H_1, \ldots, H_{\theta_n}$ are the hyperplanes of $PG(n,q)$.
Let $S$ be a subspace of $PG(n,q)$ of dimension at least $2$ and let $\mathcal{C}(S)$ be the linear code of points and hyperplanes of $S$. Then the restriction of $c$ to $S$ is defined as
$$c_{|S}=\alpha_1v^{H_1\cap S}+\cdots+\alpha_{\theta_n}v^{H_{\theta_n}\cap S}.$$
Note that $S \cap H_i$ is either a hyperplane of $S$ or is equal to $S$, for each $i$, so
$c_{|S}\in \mathcal{C}(S)$, since $c_{|S}$ is the sum of a linear combination of incidence vectors of hyperplanes of $S$ and of an $\mathbb{F}_p-$proportional of the all-one
vector $\mathbf{j}$ belonging to $\mathcal{C}(S)$. Also, $supp(c_{|S})=supp(c)\cap S$.
Furthermore, if $c=v^X$ for some subset $X\subseteq PG(n,q)$, then its restriction to a subspace $S$ is the incidence vector of $X\cap S$,
that is $c_{|S}=v^{X \cap S}$.
\end{rem}

\smallskip

Using Theorem \ref{BBW} we prove the first point of Theorem \ref{us} in the planar case for any prime $p$.

\smallskip

\begin{thm}\label{nr}
There are no codewords of $\mathcal{C}(2,q)$ with weight in the interval $]q+1,2q[$, $q=p^h$, $p$ prime.
\end{thm}
\proof
Let $c \in \mathcal{C}(2,q)$ with weight in $]q+1,2q[$.
By Theorem \ref{3.3.11}, $supp(c)$ defines a minimal blocking set $B$ of the plane $PG(2,q)$, which intersects every line in $1(mod \hspace{0.1cm} p)$
points and the nonzero components of $c$ are equal to some $a\in \mathbb{F}_p^*$.
Dividing by $a$ the codeword $c$, we obtain another codeword $c'$ of $\mathcal{C}(2,q)$.
So we may assume that the nonzero components of $c$ are $1$ and hence $c=v^B$.
If $q+k+1$ is the cardinality of $B$, by 3. of Theorem \ref{3.3.11}, we have $|B| \equiv 1 \pmod p$.
Consider a point $P \in B$. Let $t$ be a tangent line through $P$ and let $Q$ be in $t \setminus B$. Since every secant line to $B$
has at least $1+p$ points of $B$, the number of tangent lines to $B$ through $Q$ is at least $\displaystyle q-\frac{k}{p}+1$.
By Theorem \ref{BBW} the points of $B$ which belong to tangent lines through $Q$ are collinear.
Therefore there exists at least one secant line $l$ to $B$ through $P$ containing at least $\displaystyle q-\frac{k}{p}+1$ points of $B$.
Since $k<q-1$, we have that
$$ |B \cap l| \geq q-\frac{k}{p}+1> q - \frac{q}{p}+\frac{1}{p}+1, $$
and hence $|B\cap l| \geq q+1$, if $q=p$ and
\begin{equation}\label{x}
|B \cap l| \geq q - \frac{q}{p} + p +1,
\end{equation}
if $q>p$. In the first case we have $l \subseteq B$ and, by the minimality of $B$, we get $B=l$, a contradiction. Hence let $q>p$.
So we have that for each point $P \in B$ there exists a line $l$ through $P$ containing at least $\displaystyle q - \frac{q}{p} + p +1$ points of $B$.
Since $B$ is not a line and cannot be contained in the union of two lines, there exist at least $3$ lines $l_1, l_2$ and $l_3$ satisfying \ref{x}.
So,
$$|B| \geq 3 \left( q -\frac{q}{p} + p -1 \right) + 3=3q-3\frac{q}{p}+3p,$$
hence
$$ |B| \geq 3q \frac{p-1}{p}+3p+1,$$
and this is not possible if $p \geq 3$.
Finally, let $p=2$ and note that there cannot exist another line, different from $l_1, l_2$ and $l_3$ intersecting $B$ in at least
$\displaystyle \frac{q}{2}+3$ points, otherwise $ \displaystyle |B| \geq 4 \cdot \frac{q}{2}+4$.
In this way, we have shown that $B=(B\cap l_1) \cup (B \cap l_2) \cup (B \cap l_3)$ and hence, since $B$ is a blocking set and $|B|<2q$, $l_i \cap l_j \in B$ for each $i,j
\in \{1,2,3\}$ with $i \neq j$.
Let $P$ be the intersection point of $l_1$ and $l_2$ and let $Q$ be a point of $l_3\cap B$ different from $l_1\cap l_3$ and $l_2\cap l_3$, then
$|PQ \cap B|=2$ and this is not possible by point 3. of Theorem \ref{3.3.11}.
\qed

\medskip

Now, we are able to prove 1. of Theorem \ref{us} in the general case.

\begin{thm}\label{nr1}
There are no codewords of $\mathcal{C}(n,q)$ with weight in the interval $]\theta_{n-1},2q^{n-1}[$, where $q=p^h$, $p$ prime.
\end{thm}
\proof
We prove the theorem by induction on $n$. The statement holds in the case $n=2$ by Theorem \ref{nr}.
Now suppose $n>2$ and the that statement holds for each $m$ less than $n$.
Let $c \in \mathcal{C}(n,q)$ with weight in $]\theta_{n-1},2q^{n-1}[$.
By Theorem \ref{3.3.11}, $B=supp(c)$ is a minimal blocking set of $PG(n,q)$ with respect to lines and $B$ meets every line in
$1 \pmod p$ points. The nonzero components of $c$ are equal to some $a\in \mathbb{F}_p^*$,
so, up to a scalar, we may assume that the nonzero components of $c$ are $1$ and hence $c=v^B$.
Now, let $P$ be a point of $B$ and suppose that there is an integer $m$ such that $1 \leq m \leq n-2$ and there exists an $m-$subspace $S_m$
such that $S_m\cap B$ is an $(m-1)-$subspace through $P$.
In this case there exists an $(m+1)-$subspace containing $S_m$ such that $S_{m+1}\cap B$ is an $m-$subspace.
Indeed, if each $(m+1)-$subspace $S_{m+1}$, which contains $S_m$, intersects $B$ in at least $2q^{m}$ points, then
$$ 2q^{n-1}> |B| \geq \frac{\theta_n - \theta_m}{q^{m+1}}(2q^m-\theta_{m-1})+\theta_{m-1}, $$
and we get
$$ (q^{n-m-1}-1)[q^{m-1}(q-2)+1]<0, $$
a contradiction for every $q$.
So there exists an $(m+1)-$subspace $S_{m+1}$ containing $S_m$ such that $0<|B \cap S_{m+1}| < 2q^{m}$. By Remark \ref{rm}, the restriction of $c$ to $S_{m+1}$
is a codeword of $\mathcal{C}(S_{m+1})=\mathcal{C}(m+1,q)$ and its support is $B \cap S_{m+1}$, i.e. $c_{|S_{m+1}}=v^{B \cap S_{m+1}}$.
So, by the induction hypothesis and by Theorem \ref{mw} we have that $B \cap S_{m+1}$ is an $m-$subspace through $P$.
Now, since $B$ is minimal, we know that for each point $P \in B$ there exists a tangent line $l$, then we can apply the previous considerations to obtain the
existence of a hyperplane $S_{n-1}$ through $P$ such that $S_{n-1}\cap B=S_{n-2}$.
Let $\mathcal{S}$ be the set of all the $(n-2)-$subspaces $S_{n-2}$ of $PG(n,q)$ for which there exists a hyperplane $\tilde{S}_{n-1}$ such that $\tilde{S}_{n-1}\cap B=S_{n-2}$.
Note that for each point of $B$ there exists an element of $\mathcal{S}$ through it, and since $|B|>\theta_{n-1}$, it is clear that $\mathcal{S}$ contains
at least two elements.
Let $S_{n-2}$ and $S_{n-2}'$ be two elements of $\mathcal{S}$, then $S_{n-2} \cap S_{n-2}'$ is either an $(n-3)-$subspace or an $(n-4)-$subspace.
In the latter case, each hyperplane through $S_{n-2}$ has to intersect $S_{n-2}'\setminus S_{n-2}$, but this is not possible for the hyperplane $\tilde{S}_{n-1}$, since $\tilde{S}\cap B=S_{n-2}$.
Now, consider $S,S' \in \mathcal{S}$ and let $\overline{S}_{n-1}=S \vee S'$ and $\overline{S}_{n-3}=S \cap S'$.
Since the intersection of two elements of $\mathcal{S}$ is always an $(n-3)-$subspace, either $B\subseteq\overline{S}_{n-1}$ and this is not possible, or $B$ is a
cone with vertex
$\overline{S}_{n-3}$.
In this case, consider a plane $\pi$ disjoint from $\overline{S}_{n-3}$ and note that each element of $\mathcal{S}$ intersects $\pi$
in a point, hence, if $x$ is the size of $\mathcal{S}$, then
$$ |B| \geq xq^{n-2}+\theta_{n-3}.$$
Since $|B|<2q^{n-1}$, we have $0<x=|\pi \cap B|<2q$.
So, by Theorem \ref{nr}, by Theorem \ref{mw} and by Remark \ref{rm} $\pi \cap B$ is a line $r$. Then $B \subseteq \langle r, \overline{S}_{n-3} \rangle=\overline{S}_{n-1}$, a contradiction.
\qed

\smallskip

Now, we characterize the codewords with weight $2q^{n-1}$ in $\mathcal{C}(n,q)$.\\
In the planar case the following holds.

\smallskip

\begin{thm}\label{3.3.6}\cite[Corollary 6.4.4]{As}
\begin{enumerate}
  \item The minimum weight of $\mathcal{C}(2,q) \cap \mathcal{C}(2,q)^\perp$ is $2q$.
  \item The codewords of $\mathcal{C}(2,q) \cap \mathcal{C}(2,q)^\perp$ with weight $2q$ are, up to scalar, the difference of incidence vectors of any two distinct lines.
\end{enumerate}
\end{thm}

\smallskip

M. Lavrauw, L. Storme and G. Van de Voorde in \cite{Lav} generalize the first point of the previous result.

\smallskip

\begin{thm}\label{3.3.7}\cite[Theorem 5]{Lav}
The minimum weight of $\mathcal{C}(n,q) \cap \mathcal{C}(n,q)^\perp$ is $2q^{n-1}$.
\end{thm}

\smallskip

Also, if $q=p$ the words in $\mathcal{C}(n,p)\cap \mathcal{C}(n,p)^\perp$ with weight $2p^{n-1}$ are the scalar multiples of the difference of the incidence vectors of two hyperplanes of $PG(n,p)$, see \cite[Remark 2]{Lav} and \cite[Theorem 12]{Lav2}.

\smallskip

In the next, we will prove that the second minimum weight of $\mathcal{C}(n,q)$ is $2q^{n-1}$, $q=p^h$, $h\geq 1$ and that the words of $\mathcal{C}(n,q)$ with this weight are, up to scalar, the difference of the incidence vectors of two hyperplanes of $PG(n,q)$ for each prime $p$.

\smallskip

In the same way as it was done by the authors in \cite[Lemma 6]{Lav}, one can prove the following.

\smallskip

\begin{prop}\label{2}
There are no codewords in $\mathcal{C}(n,q) \setminus \mathcal{C}(n,q)^\perp$ with weight $2q^{n-1}$.
\end{prop}
\proof
Suppose to the contrary that there exists $c$ be in $\mathcal{C}(n,q) \setminus \mathcal{C}(n,q)^\perp$ with weight $2q^{n-1}$.
Since $c \notin \mathcal{C}(n,q)^\perp$, by 2. and 3. of Property \ref{3.3.5}, for each subspace $U$ with $dim \hspace{0.1cm} U \geq 1$ we have $(c,v^U) = a$,
for some $a \in \mathbb{F}_p^*$, i.e. $U \cap supp(c) \neq \emptyset$.
In particular, this holds for the lines, and so $B=supp(c)$ is an $(n-1)-$blocking set in $PG(n,q)$.
Also, if $R$ is an essential point of $B$ and $t$ is a tangent line to $B$ through $R$, then $a$ is the component of $c$ corresponding to $R$.
This implies that $B$ is minimal, indeed if there exists a point $P \in B$ which is not an essential point for $B$, since $|B|=2q^{n-1}$,
there exists a line $l$ through $P$
intersecting $B$ in exactly two points. If $l\cap B=\{P,Q\}$, then by Theorem \ref{min}, $Q$ is an essential point of $B$.
So, the corresponding component of $c$ is $a$ and, denoted by $x$
the component of $c$ corresponding to $P$, we get
$$ (c,v^{l})=x+a = a,$$
i.e. $x = 0$, a contradiction.
So $B$ is a minimal $(n-1)-$blocking set and hence for each point of $B$ there exists a tangent line.
This means that the nonzero components of $c$ are equal to $a$, that is $c=av^B$.
Since $c \in \mathcal{C}(n,q)\setminus \mathcal{C}(n,q)^\perp$, we have that $v^B \in \mathcal{C}(n,q)\setminus \mathcal{C}(n,q)^\perp$.
Then, each line of $PG(n,q)$ intersects $B$ in $1(mod \hspace{0.1cm} p)$ points, and hence
$|B| \equiv 1 \pmod p$, but this is not possible since $|B|=2q^{n-1}$.
\qed

\smallskip

\begin{lem}\label{lem}
Let $X \subseteq PG(n,q)$, $n \geq 2$, with $|X|=2q^{n-1}$ and such that for each $h-$dimensional subspace $S_h$, with $1 \leq h \leq n-1$, one of the following occurs:
\begin{enumerate}
  \item $X\cap S_h = \emptyset$;
  \item $X\cap S_h$ is the symmetric difference of two hyperplanes of $S_h$ (if $h=1$, $S_1$ is a $2-$secant line to $X$);
  \item $X\cap S_h=S_h\setminus S_{h-1}$, where $S_{h-1}$ is a hyperplane of $S_h$;
\end{enumerate}
then $X$ is the symmetric difference of two hyperplanes of $PG(n,q)$.
\end{lem}
\proof
Note that if $q=2$, then points 2. and 3. describe the same set of points and hence for each hyperplane $S_{n-1}$ of $PG(n,2)$ we have that either $S_{n-1}\cap X=\emptyset$ or $S_{n-1}\setminus S_{n-2}\subseteq X$ where $S_{n-2}$ is a hyperplane of $S_{n-1}$.
Since $|X|=2q^{n-1}=2^n$ we easily get that $X$ is the symmetric difference of two hyperplanes of $PG(n,2)$.
So, let $q>2$ and we prove that:
\begin{itemize}
  \item[(a)] For each $P \in X$ there exists a line $l$ through $P$ such that $l \setminus \{Q\}\subseteq X$, for some $Q \in l$;
  \item[(b)] If $S_m$ is an $m-$subspace of $PG(n,q)$, with $0 \leq m \leq n-2$, for which case 3. holds, then there exists an $(m+1)-$subspace $S_{m+1}$ containing $S_m$ satisfying 3..
\end{itemize}
Let $P$ be a point of $X$ and assume that (a) is not satisfied, then every line through $P$ is a $2-$secant line to $X$. Hence,
$$ |X|=\theta_{n-1}+1, $$
and this is not possible since $q>2$.
So, (a) is proved.
Now, let $S_m$ be an $m-$subspace, with $0 \leq m \leq n-2$, such that $S_m \setminus S_{m-1}\subseteq X$, where $S_{m-1}$ is a hyperplane of $S_m$,
and assume that every $S_{m+1}$ containing $S_m$ intersects $X$ in the symmetric difference of two $m-$subspaces, one of which is $S_m$, then
$$|X|=\theta_{n-m-1}(2q^m-q^m)+q^m,$$
and since $m \leq n-2$ this is possible only for $q=2$.
In this way we have proven (b).
Now, from (a) and (b) we get that for each point $P \in X$
there exists an $S_{n-1}$ through $P$ such that
$S_{n-1}\setminus S_{n-2}\subseteq X$.
If there exists another hyperplane $\overline{S}_{n-1}$ containing $S_{n-2}$ such that $\overline{S}_{n-1}\setminus S_{n-2} \subseteq X$, then $X$ is the symmetric difference of two hyperplanes.
Otherwise, denoted by
$x$ the number of the hyperplanes through $S_{n-2}$ intersecting $X$ in a symmetric difference of two of its hyperplanes, we get
$$q^{n-1}=x2q^{n-2},$$
and so $\displaystyle x=\frac{q}{2}$, and this is not possible if $p\geq 3$.
Finally, let $p=2$ and $q>2$, i.e. $q \geq 4$.
If there exist at least $3$ hyperplanes satisfying 3., we obtain
$$ |X| \geq q^{n-1}+q^{n-1}-q^{n-2}+q^{n-1}-2q^{n-2},$$
which is not possible for $q > 3$.
Then there exist two hyperplanes verifying 3., and since $|X|=2q^{n-1}$, $X$ is the symmetric difference of two hyperplanes of $PG(n,q)$.
\qed

\smallskip

\begin{thm}
The codewords in $\mathcal{C}(n,q)$ of weight $2q^{n-1}$ are the scalar multiples of the difference of the incidence vectors of two distinct
hyperplanes of $PG(n,q)$.
\end{thm}
\proof
The assert holds in the case $n=2$ by Theorem \ref{3.3.6} and by Proposition \ref{2}.
Now, suppose the assert true in the code $\mathcal{C}(t,q)$, with $2 \leq t \leq n-1$ and
let $c \in \mathcal{C}(n,q)$ with $wt(c)=2q^{n-1}$. Then, by Proposition \ref{2}, $c\in\mathcal{C}(n,q)\cap\mathcal{C}(n,q)^\perp$.
Denote by $X$ the support of $c$ and note that, by 3. of Property \ref{3.3.3}
$(c,v^U)\equiv 0 \pmod p$ for each subspace $U$ of $PG(n,q)$ of dimension $h$, with $1 \leq h \leq n-1$, hence $X$ has no tangent space.
Also, if $U$ is a subspace of dimension $h$ with $2\leq h \leq n-1$, then by Remark \ref{rm} $supp(c_{|U})=supp(c)\cap U$ and $c_{|U}\in \mathcal{C}(U)
=\mathcal{C}(h,q)$, hence if $supp(c)\cap U\neq \emptyset$ by the Theorem \ref{nr1} and by the induction hypothesis one of the following holds:
\begin{itemize}
  \item[(a)] $supp(c) \cap U$ is a hyperplane of $U$;
  \item[(b)] $supp(c) \cap U$ is a symmetric difference of two hyperplanes of $U$;
  \item[(c)] $|supp(c)\cap U|>2q^{h-1}$.
\end{itemize}
Now, we are able to prove the following:
\begin{itemize}
  \item[$(\ast)$] If $S_m$ is an $m-$subspace, with $1 \leq m \leq n-2$, such that $S_m \cap X$ is the symmetric difference of two hyperplanes of $S_m$ (if $m=1$ $S_1$ is a $2-$secant line to $X$), then each $(m+1)-$subspace $S_{m+1}$
containing $S_m$ intersects $X$ in the symmetric difference
of two $m-$subspaces.
\end{itemize}
Indeed,
if $S_{m+1}$ is an $(m+1)-$subspace containing $S_m$, since $S_{m+1}\cap X \supseteq S_m \cap X$, Case (a) does not occur, hence
either $S_{m+1}\cap X$ is the symmetric difference of two hyperplanes of $S_{m+1}$ or $|S_{m+1}\cap X|>2q^m$.
Let $x$ be the number of $S_{m+1}$
containing $S_m$ such that $ S_{m+1}\cap X$ is the symmetric difference of two hyperplanes of $S_{m+1}$, then:
$$ 2q^{n-1}=|X|\geq x(2q^{m}-2q^{m-1})+\left(\theta_{n-m-1}-x\right)(2q^m-2q^{m-1}+1)+2q^{m-1}, $$
where $\theta_{n-m-1}$ is the number of the $(m+1)-$subspaces containing $S_m$. Then we get
$$ x \geq \theta_{n-m-1},$$
and hence $x=\theta_{n-m-1}$, i.e. each $(m+1)-$subspace $S_{m+1}$
containing $S_m$ intersects $X$ in the symmetric difference
of two $m-$subspaces.
Since $|X|=2q^{n-1}$ and there are no tangent lines to $X$, for each point $P$ of $X$ there exists a $2-$secant line, so applying $(\ast)$,
we get that for each point $P \in X$, each $h-$subspace, with $2 \leq h \leq n-1$,
through a $2-$secant line containing $P$ intersects $X$ in the symmetric difference of two $(h-1)-$subspaces,
one of which contains $P$.
As a consequence, we get that for a line $l$ of $PG(n,q)$ one of the following holds true:
\begin{itemize}
  \item[(a')] $l$ is external to $X$;
  \item[(b')] $l$ is $2-$secant to $X$;
  \item[(c')] $l\setminus \{Q\}$ is contained in $X$, where $Q \in l$.
\end{itemize}
Indeed, if $l$ is a line which contains more than two points and $R\in l \cap X$, through $R$ there exists at least one
$2-$secant line $l'$ to $X$. The plane $l \vee l'$ intersects $X$ in the symmetric difference of two lines, and one of these must be the line $l$.\\
By the previous considerations, moving forward by finite induction on $h$, we get that for each subspace $S_h$ of $PG(n,q)$ with $1 \leq h \leq n-1$ one of the following occurs:
\begin{enumerate}
  \item $S_h$ is external to $X$;
  \item $X\cap S_h$ is the symmetric difference of two hyperplanes of $S_h$ (if $h=1$, $S_1$ is a $2-$secant line to $X$);
  \item $S_h\setminus S_{h-1}\subseteq X$, where $S_{h-1}\subseteq S_h$;
\end{enumerate}
then, by Lemma \ref{lem}, $X$ is the symmetric difference of two distinct hyperplanes.
\qed

\end{document}